\documentclass[a4paper]{article}
\usepackage{url}

\pdfoutput=1
\author{Robin Houston}
\title{Tackling the Minimal Superpermutation Problem}

\begin{document}
\maketitle

\begin{abstract}
\noindent A superpermutation on $n$ symbols is a string that contains each of the $n!$ permutations of the $n$ symbols as a contiguous substring. The shortest superpermutation on $n$ symbols was conjectured to have length $\sum_{i=1}^n i!$. The conjecture had been verified for $n \leq 5$. We disprove it by exhibiting an explicit counterexample for $n=6$. This counterexample was found by encoding the problem as an instance of the (asymmetric) Traveling Salesman Problem, and searching for a solution using a powerful heuristic solver.
\end{abstract}

\section{The Minimal Superpermutation Problem}\label{section:msp}

A superpermutation on $n$ symbols is a string that contains each of the $n!$ permutations of the $n$ symbols as a contiguous substring. For the sake of concreteness, we'll assume that the $n$ symbols are the numbers $1, 2, \dots, n$. The problem of finding minimal-length superpermutations was posed by Ashlock and Tillotson \cite{AshlockTillotson}, who established that $\sum_{i=1}^n i!$ symbols suffice and conjectured that this is optimal.

There is a simple construction that takes a superpermutation on $n-1$ symbols of length $k$ and produces a superpermutation on $n$ symbols of length $k + n!$. Iterating this construction gives a sequence of palindromic superpermutations of length $1!, 1!+2!, 1!+2!+3!, \dots$, viz:

\begin{itemize}
\item \texttt{1}
\item \texttt{121}
\item \texttt{123121321}
\item \texttt{123412314231243121342132413214321}
\item \texttt{1234512341\allowbreak 5234125341\allowbreak 2354123145\allowbreak 2314253142\allowbreak 3514231542\allowbreak 3124531243\allowbreak 5124315243\allowbreak 1254312134\allowbreak 5213425134\allowbreak 2153421354\allowbreak 2132451324\allowbreak 1532413524\allowbreak 1325413214\allowbreak 5321435214\allowbreak 3251432154\allowbreak 321}
\item \dots
\end{itemize}

\noindent The construction is the following. Order the $(n-1)!$ permutations of $n-1$ symbols according to the order they first appear in the superpermutation; then replace each permutation $s$ by the $n$ permutations of $n$ symbols that can be obtained as cyclic shifts of $sn$, shifting leftwards one position at a time from $sn$ to $ns$. Finally eliminate overlap between adjacent permutations to obtain a superpermutation.

For example, the superpermutation \texttt{121} on two symbols gives the sequence \texttt{12}, \texttt{21}. The permutation \texttt{12} is expanded to \texttt{123}, \texttt{231}, \texttt{312}, and the permutation \texttt{21} expanded to \texttt{213}, \texttt{132}, \texttt{321}, so the final sequence of permutations is \texttt{123}, \texttt{231}, \texttt{312}, \texttt{213}, \texttt{132}, \texttt{321}. Overlaps are eliminated to yield \texttt{123121321}.

It is easy to show, by exhaustive enumeration using a computer, that for $n\leq 4$ these palindromic superpermutations have minimal length and are unique of this length up to relabelling. Johnston \cite{Johnston} showed that uniqueness fails for $n\geq 5$: there is more than one superpermutation of length $\sum_{i=1}^n i!$. More recently Benjamin Chaffin used a clever exhaustive computer search \cite{OEIS:A180632, JohnstonMinimal5} to find all eight five-symbol superpermutations of length $\sum_{i=1}^5 i!=153$, and show that there are no shorter ones.

The minimal length is still unknown for $n\geq 6$, but we can show that for all $n\geq 6$ it is strictly less than the conjectured length $\sum_{i=1}^n i!$.

\section{The Travelling Salesman Problem}\label{section:tsp}

The Travelling Salesman Problem (TSP) is the problem of finding a minimum-weight Hamiltonian circuit of a weighted graph. We are interested mainly in the asymmetric TSP, where the graph is directed.

For an $n$-symbol alphabet, we may construct a complete directed graph with $n!$ vertices, one for each permutation of the $n$ symbols. The weight of the edge from $s$ to $t$ is the least $0\leq k\leq n$ such that the $(n-k)$-suffix of $s$ is equal to the $(n-k)$-prefix of $t$. Then a minimum-weight Hamiltonian path in this graph corresponds to a minimal superpermutation on $n$ symbols. The length of the corresponding superpermutation is $n+w$ where $w$ is the weight of the minimum-weight Hamiltonian path.

So far this is not quite an instance of the TSP, since we are looking for a Hamiltonian \emph{path} rather than a Hamiltonian circuit. To relate it to the TSP it is sufficient to change the weights of certain edges to $0$. Specifically, let $o$ be the identity permutation $12\dots n$ and let the weight of each edge $s$--$o$ be $0$. A minimum-weight Hamiltonian circuit on the resulting graph corresponds to a minimum-weight Hamiltonian path starting at $o$.

The TSP has been intensively studied, and though it is NP-hard -- hence no worst-case polynomial-time algorithm is known -- there are solvers that in practice work very well. In particular we used Concorde \cite{Concorde}, which finds provably-minimal solutions; and LKH \cite{LKH2000, LKH2009}, which is a fast, randomised, approximate solver.

Concorde only works for the symmetric TSP, so to apply Concorde it was necessary to use the Jonker-Volgenant construction \cite{JonkerVolgenant, JonkerVolgenantErratum} to convert our instance of the asymmetric TSP on $n!$ vertices to an instance of the symmetric TSP on $2(n!)$ vertices.

We were able to use Concorde to solve the five-symbol instance in 1929.02 seconds on an Amazon EC2 `m3.medium' instance running Linux, confirming Chaffin's recent result. We have not found an exact solution to the six-symbol instance -- Concorde failed with an internal error after running for several days. Nevertheless we conjecture it is within reach of current technology.

On the other hand, we were able to falsify the minimal-length conjecture on six symbols without solving the instance completely, using LKH to search repeatedly for approximate solutions till we found a Hamiltonian circuit of weight 866, representing a superpermutation of length 872 (less than the conjectured minimum of 873). Since LKH supports the asymmetric TSP directly, it was not necessary to use the larger symmetrical instance in this case. The first superpermutation of length 872 was found on the 9228th run -- with each run taking approximately 2--4 seconds on a 2010 MacBook Pro with 2.66 GHz Intel Core 2 Duo processor. The parameters used were:\\[1em]\noindent
\texttt{BACKTRACKING = YES\\
MAX\_CANDIDATES = 6 SYMMETRIC\\
MOVE\_TYPE = 3\\
PATCHING\_C = 3\\
PATCHING\_A = 2}
\\[1em]\noindent These parameters are those used by Helsgaun to find optimal solutions to David Soler's ATSP instances \cite{HelsgaunSoler}, with the addition of the \texttt{BACKTRACKING} option.

\section{A short superpermutation on six symbols}\label{section:short}

The following superpermutation on six symbols has length 872, which is less than the conjectured minimum of
$\sum_{i=1}^6 i!=873$. The recursive construction described in Section~\ref{section:msp} may then be used to construct a superpermutation shorter than the conjectured minimum for all $n > 6$. \\

\noindent\texttt{1234561234\allowbreak 5162345126\allowbreak 3451236451\allowbreak 3264513624\allowbreak 5136425136\allowbreak 4521364512\allowbreak 3465123415\allowbreak 6234152634\allowbreak 1523641523\allowbreak 4615234165\allowbreak 2341256341\allowbreak 2536412534\allowbreak 6125341625\allowbreak 3412653412\allowbreak 3564123546\allowbreak 1235416235\allowbreak 4126354123\allowbreak 6541326543\allowbreak 1264531624\allowbreak 3516243156\allowbreak 2431652431\allowbreak 6254316245\allowbreak 3164253146\allowbreak 2531426531\allowbreak 4256314253\allowbreak 6142531645\allowbreak 2314652314\allowbreak 5623145263\allowbreak 1452361452\allowbreak 3164532164\allowbreak 5312643512\allowbreak 6431526431\allowbreak 2564321564\allowbreak 2315462315\allowbreak 4263154236\allowbreak 1542316542\allowbreak 3156421356\allowbreak 4215362415\allowbreak 3621453621\allowbreak 5436215346\allowbreak 2135462134\allowbreak 5621346521\allowbreak 3462513462\allowbreak 1536421563\allowbreak 4216534216\allowbreak 3542163452\allowbreak 1634251634\allowbreak 2156432516\allowbreak 4325614325\allowbreak 6413256431\allowbreak 2654321654\allowbreak 3261534261\allowbreak 3542613452\allowbreak 6134256134\allowbreak 2651342615\allowbreak 3246513246\allowbreak 5312463512\allowbreak 4631524631\allowbreak 2546321546\allowbreak 3251463254\allowbreak 1632546132\allowbreak 5463124563\allowbreak 2145632415\allowbreak 6324516324\allowbreak 5613245631\allowbreak 2465321465\allowbreak 3241653246\allowbreak 1532641532\allowbreak 6145326154\allowbreak 3265143625\allowbreak 1436521435\allowbreak 6214352614\allowbreak 3521643521\allowbreak 4635214365\allowbreak 1243615243\allowbreak 6125436124\allowbreak 5361243561\allowbreak 2436514235\allowbreak 6142351642\allowbreak 3514623514\allowbreak 2635142365\allowbreak 1432654136\allowbreak 2541365241\allowbreak 3562413526\allowbreak 4135246135\allowbreak 2416352413\allowbreak 6542136541\allowbreak 23}

\section{Remarks}
\noindent General-purpose solvers for NP-hard problems -- SAT, SMT, TSP, etc. -- are remarkably powerful and applicable to many combinatorial problems. The world of pure mathematics has been slow to exploit this opportunity, but that is beginning to change: the recent advances on the Erd\H{o}s Discrepancy Problem using a SAT solver \cite{konev2014, konev2014sat} are especially notable in this regard.

It seems plausible that the minimal superpermutation problem for $n=6$ could be solved exactly using Concorde or a similar algorithm with a few weeks or months of CPU time. More generally, these instances seem to present an interesting challenge for the TSP community.

\section{Ancillary files}
Several ancillary files are included, listed with the arXiv abstract.

\subsection*{Python scripts and their output}
\begin{itemize}
	\item \texttt{mkatsp.py} Generate an instance of the asymmetric TSP representing the
	minimal superpermutation problem for a given value of $n$. The output is a TSPLIB file
	of type ATSP. The files \texttt{5.atsp} and \texttt{6.atsp} are the output of this script
	for $n=5$ and $n=6$.

	\item \texttt{symmetrise.py} Apply the Jonker-Volgenant transformation to transform
	an ATSP instance into a corresponding symmetric TSP instance. The files \texttt{5.tsp}
	and \texttt{6.tsp} are the result of applying this to \texttt{5.atsp} and \texttt{6.atsp}
	respectively.

	\item \texttt{mkpalindromic.py} Produce a simple palindromic superpermutation on the
	specified number of symbols.

	\item \texttt{printtour.py} Take a tour produced by LKH and print the corresponding superpermutation.

	\item \texttt{splitsuperperm.py} Split a superpermutation into permutations. This is useful for analysing the structure of a superpermutation, and (by sorting the output) for verifying that a supposed superpermutation really does produce all permutations.
\end{itemize}

\subsection*{Other files}
\begin{itemize}
	\item \texttt{concorde-output-5-symbols.txt} The log of running Concorde on \texttt{5.tsp}.
	\item \texttt{6.par} The parameter file used to run LKH on \texttt{5.atsp}.
	\item \texttt{6.866.lkh} The first tour of weight 866 obtained when running LKH with parameter file \texttt{5.par}.
	\item \texttt{superperm-6-866.txt} The superpermutation of length 872, obtained by running \texttt{printtour.py} on \texttt{6.866.lkh}.
\end{itemize}

\bibliographystyle{plain}
\bibliography{refs}
\end{document}